\documentclass[prb,article,nofootinbib,11pt]{revtex4-2} 

\usepackage{amsmath}
\usepackage{amssymb}
\usepackage{amsfonts} 
\usepackage{graphicx} 
\usepackage{physics}
\usepackage{subfigure}
\usepackage{mathrsfs}
\usepackage{accents}
\usepackage{cancel}
\usepackage{longtable}
\usepackage[usenames,dvipsnames]{xcolor}
\usepackage{mathtools}
\usepackage{relsize}
\usepackage{enumerate}
\usepackage{xfrac}
\usepackage{bm}
\usepackage{nicefrac}
\usepackage{listings}

\begin{document}


\title{On the Equilibria and Bifurcations of a Rotating Double Pendulum}

\author{Jonathan Tot}
\email{jonathan.tot@dal.ca} 
\affiliation{Department of Mathematics and Statistics, Dalhousie University, Halifax, NS B3H-4R2}

\author{Robert H. Lewis}
\email{rlewis@fordham.edu}
\affiliation{Fordham University, Bronx, NY 10458}


\begin{abstract}
	The double pendulum, a simple system of classical mechanics, is widely studied as an example of, and testbed for, chaotic dynamics.  In \cite{maiti2016}, Maiti et al. study a generalization of the simple double pendulum with equal point-masses at equal lengths, to a rotating double pendulum, fixed to a coordinate system uniformly rotating about the vertical.  In this paper, we study a considerable generalization of the double pendulum, constructed from physical pendula, and ask what equilibrium configurations exist for the system across a comparatively large parameters space, as well as what bifurcations occur in those equilibria.  Elimination algorithms are employed to reduce systems of polynomial equations, which allows for equilibria to be visualized, and also to demonstrate which models within the parameter space exhibit bifurcation.  We find the \texttt{DixonEDF} algorithm for the Dixon resultant~\cite{lewis2017}, written in the computer algebra system (CAS) \textit{Fermat} \cite{fermat}, to be capable to complete the computation for the challenging system of equations that represents bifurcation, while attempts with other algorithms were terminated after several hours.
\end{abstract}

\maketitle

\section{Introduction} 

This work is concerning a system of classical physics, namely a rotating double pendulum (RDP). Here, ‘rotating’ is used to indicate a double pendulum that is made to rotate about a vertical axis with a constant angular velocity $\omega_a$. A double pendulum has two joints or pivots, and we take the vertical axis of this rotation to pass through the inner, stationary pivot. Additionally, our consideration is of a physical double pendulum, constructed from two 3-dimensional rigid bodies. Thus the physics depends on many dimensional parameters: masses, lengths, $\omega_a$, the strength of gravity, and principal moments of inertia – 13 in total. However, only 6 dimensionless parameters are required to describe the dynamics of the system.

\subsection{Double Pendulum Systems}

Many researchers have studied similar systems. Studies of bifurcations in the non-linear dynamics of double pendula can be found in \cite{samaranyake1993,yu_bi1998,qinsheng2000}. In 2001, Bridges and Georgiou \cite{bridges2001} studied a transverse rotating double pendulum, in which the axes of the two pivots are not parallel, so that the pendula do not swing in the same plane.  This system admits a doubly-degenerate equilibrium at the trivial solution (both pendula hanging down) in linearization, indicated by the coalescence of four zero eigenvalues. It is claimed that this is the simplest autonomous system one could construct with two degrees of freedom and admits such a critical point. In \cite{maiti2016}, S. Maiti et. al. study a rotating double pendulum model of equal point-masses on equal length massless rods. With increasing rotation speed, they observe a transition from chaotic dynamics, to quasi-periodic order, and back to chaos, as evidenced by Poincaré sections.

\subsection{Structure of this paper}

In \S~\ref{S1} we introduce the physical construction of the rotating double pendulum that we consider, and introduce the system's Lagrangian (which is derived in Appendix~\ref{Ap.RDP.deriv}) as well as the dimensionless parametrization that will be used, and a particular special case of point-masses on massless rods (PMMR).  In \S~\ref{eqns} we analyze the Euler-Lagrange equations to derive equations for the equilibrium configurations of the system and for bifurcation of those equilibria.  In \S~\ref{trivS} we consider certain trivial equilibrium solutions (among other respects, they are equilibria of the system for all parameter values) and analyze what bifurcations occur for these constant equilibria.  Finally, in \S~\ref{elim} we use polynomial elimination algorithms, including \texttt{DixonEDF} \cite{lewis2017} for the Dixon resultant, to  1) visualize bifurcation diagrams in the PMMR special case and 2) derive a polynomial condition on the parameter space for bifurcation.  This not only recovers the bifurcations for the trivial equilibria, but confirms and describes an addition non-trivial bifurcation that is seen in the bifurcations plots (Fig.~\ref{bifresCP}). We then present some conclusions and some possible directions for future study.  The appendices include: the normal mode analysis of a single rotating (physical) pendulum in Appendix \ref{rpp}.  This brief analysis is very much the motivation and model for the desired analysis of the RDP---if possible, one would hope to be able to complete a normal mode analysis for the equilibrium solutions of the RDP.  Appendix \ref{Ap.RDP.deriv} details the derivation of the RDP Lagrangian, while Appendix \ref{outputs} contains long-form output of the polynomial systems for equilibrium and bifurcation that derive from this work.

\section{The Rotating Double Pendulum (RDP)}\label{S1}

\begin{figure}[!ht]
	\centering
	\includegraphics[scale=.25]{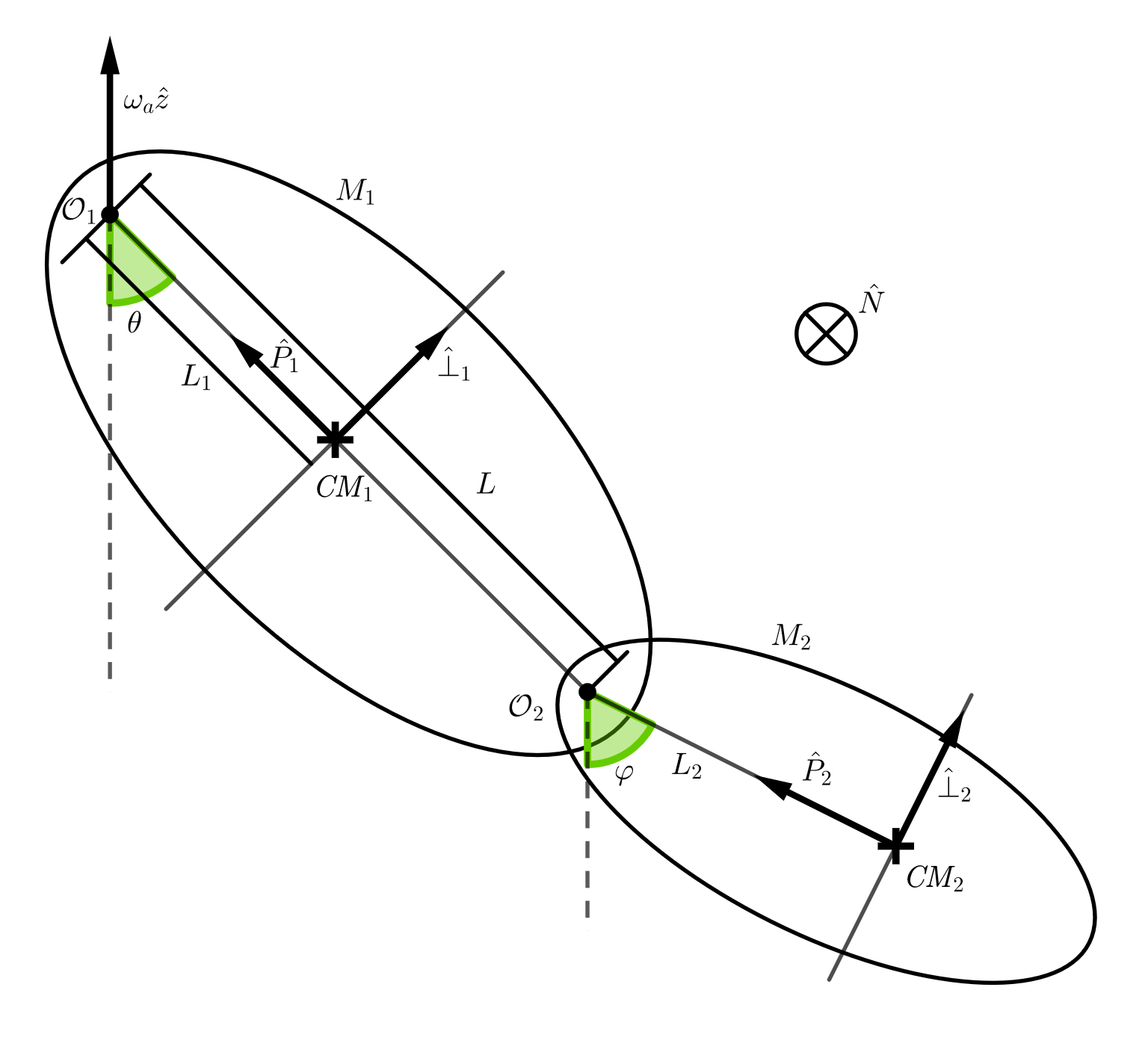}
	\caption{Schematic showing the construction and configuration of a rotating double pendulum from two physical pendula.}
	\label{RDPconfig1}
\end{figure}

We study a rotating double pendulum consisting of two physical pendula with masses $M_1$ and $M_2$.  Each pendulum has principal moments of inertia, $I_{P,\perp,N}^{(1),(2)}$, and the pendula are aligned by their principal axes, as demonstrated in Fig.~\ref{RDPconfig1}.  In particular, the pivots $\mathcal{O}_1$ and $\mathcal{O}_2$ are co-linear with $\mathit{CM}_{\!1}$ the center-of-mass of $M_1$, lying along the $P_1$ axis.  Similarly, the second pivot $\mathcal{O}_2$ and center-of-mass $\mathit{CM}_{\!2}$ of $M_2$ lie along the $P_2$ axis.  The upper pivot $\mathcal{O}_1$ is fixed, while the lower pivot $\mathcal{O}_2$ fastens the pendula together.  The $P_1$ and $P_2$ axes make angles $\theta$ and $\varphi$ with the vertical, respectively. These axes also define the plane of the pendulum -- the vertical plane in which the masses swing independently -- and this plane is made to rotate about the vertical axis through $\mathcal{O}_1$ with angular frequency $\omega_a$ (this rotation is taken to be counter-clockwise when the pendulum is viewed from above).  Further principle axes $\perp_1$, $\perp_2$ are also in the vertical plane, so that the third principle axes of the pendula coincide, normal to the plane of the pendulum (pointing into the page) $N_1=N_2=N$.

This system has the Lagrangian\footnote{The subscript 2 in $\mathcal{L}_2$ is to indicate the Lagrangian for a \textit{double} pendulum. At (\ref{rpp.lag}), $\mathcal{L}_1$ is used for the Lagrangian of a (single) rotating physical pendulum.}

\begin{align}
	\mathcal{L}_2=&\frac{1}{2}\left(A\dot{\theta}^2+2B\cos(\varphi-\theta)\dot{\theta}\dot{\varphi}+C\dot{\varphi}^2 \right) +\frac{1}{2}\left(\bar{A}\sin^2\theta+2\bar{B}\sin\theta\sin\varphi+\bar{C}\sin^2\varphi\right) \label{rdp.lag} \\
	&+K_1\cos\theta+K_2\cos\varphi \nonumber
\end{align}

\noindent where dot-derivatives indicate differentiation with respect to dimensionless time $\tau=\omega t$, the characteristic frequency $\omega$ to be determined.  The coefficients of the Lagrangian (\ref{rdp.lag}) are

\begin{align}
	A&=\left(M_1 L_1^2+M_2L^2+I^{(1)}_N\right)\omega^2 &  \bar{A}&=\left(M_1L_1^2+M_2L^2+I^{(1)}_\perp-I^{(1)}_P\right)\omega_a^2 \nonumber \\
	B&=M_2 L_2 L \omega^2 &  \bar{B}&=M_2 L_2 L \omega_a^2 \label{coeffs} \\
	C&=\left(M_2 L_2^2+I^{(2)}_N\right)\omega^2 &   \bar{C}&=\left(M_2 L_2^2+I^{(2)}_\perp-I^{(2)}_P\right)\omega_a^2 \nonumber \\
	K_1&=\left(M_1L_1+M_2L\right)g & K_2&=M_2L_2g\,. \nonumber
\end{align}

\noindent The coefficients all have dimensions of energy.  The terms of the first parentheses of (\ref{rdp.lag}), quadratics in time-derivatives, are the kinetic terms, and we refer to $A,B,C$ as the kinetic coefficients.  The terms of the second parentheses give rise to centrifugal forces; we call $\bar{A},\bar{B},\bar{C}$ then centrifugal coefficients.  The second line of (\ref{rdp.lag}) is the negative of the total gravitational potential energy of the pendula, so we will call these the gravitational terms, and $K_1,K_2$ the gravitational coefficients.

We determine $\omega$ by reducing to a single rotating pendulum: if we lock the second pivot such that $\varphi\equiv\theta$, we have

\begin{equation}
	\mathcal{L}_1=\mathcal{L}_2\vert_{\varphi\rightarrow\theta}=\frac{1}{2}\left(A+2B+C\right)\dot{\theta}^2+\frac{1}{2}\left(\bar{A}+2\bar{B}+\bar{C}\right)\sin^2\theta+\left(K_1+K_2\right)\cos\theta \label{rpp.lag}
\end{equation}

\noindent and we make the identifications

\begin{align}
	E&=A+2B+C=\left[M_1 L_1^2+M_2\left(L+L_2\right)^2+I^{(1)}_N+I^{(2)}_N\right]\omega^2 \\
	\bar{E}&=\bar{A}+2\bar{B}+\bar{C}=\left[M_1 L_1^2+M_2\left(L+L_2\right)^2+I^{(1)}_\perp+I^{(2)}_\perp-I^{(1)}_P-I^{(2)}_P\right]\omega_a^2 \\
	K&=K_1+K_2=\left[M_1 L_1+M_2(L+L_2)\right]g\,.
\end{align}

\noindent We now fix $\omega$ by the condition $E=\bar{E}+K$, which gives

\begin{equation}
	\omega^2=\frac{\bar{E}+K}{M_1 L_1^2+M_2\left(L+L_2\right)^2+I^{(1)}_N+I^{(2)}_N},
\end{equation}

\noindent and define the bounded dimensionless parameter

\begin{equation}
	Q=\frac{\bar{E}}{E}=\frac{\bar{E}}{\bar{E}+K}\in[0,1]\,. \label{defQ}
\end{equation}

\noindent Then if we normalize the Lagrangian (\ref{rpp.lag}), dividing through by $E$, we have

\begin{equation}
	\mathbf{L}_1=\mathcal{L}_1/E=\frac{1}{2}\dot{\theta}^2+\frac{1}{2}Q\sin^2\theta+(1-Q)\cos\theta\,.
\end{equation}

\noindent The equilibrium, bifurcation, and normal mode analyses of this sub-system are elementary and contained in Appendix \ref{rpp}; we propose to undertake the corresponding analysis for the RDP, to such extent as is tractable.  Next, we continue to parametrize the Lagrangian (\ref{rdp.lag}) with dimensionless combinations of the coefficients.

\subsection{Dimensionless Parameters for the Rotating Double Pendulum}

We define the following dimensionless parameters in terms of the coefficients (\ref{coeffs})

\begin{alignat}{2}
	\tilde{\delta}&=\frac{A+C-2B}{A+2B+C}\quad  &&\delta=\frac{\bar{A}+\bar{C}-2\bar{B}}{\bar{A}+2\bar{B}+\bar{C}}  \label{delta}\\
	\tilde{\sigma}&=\frac{A-C}{A+C}  &&\sigma=\frac{\bar{A}-\bar{C}}{\bar{A}+\bar{C}}  \label{sigma}
\end{alignat}
\vspace{-14.5pt}
\begin{align}
	\alpha&=\frac{(A+C)\bar{B}-(\bar{A}+\bar{C})B}{(A+2B+C)\bar{B}+(\bar{A}+2\bar{B}+\bar{C})B} \\
	\eta&=\frac{(A-C)\bar{B}-(\bar{A}-\bar{C})B}{(A+2B+C)\bar{B}+(\bar{A}+2\bar{B}+\bar{C})B} \label{eta} \\
	\chi&=\frac{K_1-K_2}{K_1+K_2} \label{chi}
\end{align}

\noindent Each of the parameters (\ref{delta}-\ref{chi}) takes values in the interval $[-1,1]$.  These 7 rational functions of the coefficients are not all independent; $\tilde{\delta}$ and $\tilde{\sigma}$ can be written in terms of $\delta,\sigma,\alpha$ and $\eta$.  We now reduce the Lagrangian (\ref{rdp.lag}), dividing by $E/4$ for convenience, and we have 

\begin{align}
	\mathbf{L}_2=&4\mathcal{L}_2/E \nonumber \\
	=&\frac{1}{2}\left[\tilde{a}\,\dot{\theta}^2+2\tilde{b}\,\cos{(\varphi-\theta)}\dot{\theta}\dot{\varphi}+\tilde{c}\,\dot{\varphi}^2\right]+\frac{Q}{2}\left(a\sin^2\theta+2b\sin\theta\sin\varphi+c\sin^2{\varphi}\right) \\
	&+2\,(1-Q)\left[(1+\chi)\cos\theta+(1-\chi)\cos\varphi\right] \nonumber
\end{align}

\noindent where $Q$ is as defined in (\ref{defQ}), and the reduced Lagrangian coefficients\footnote{For easier computations, the expressions for $\tilde{a},\tilde{b},\tilde{c}$ can be somewhat simplified by using alternate functions of $\alpha$ and $\eta$: $x=(\alpha+\eta)/(1+\alpha)$, and $y=(\alpha-\eta)/(1+\alpha)$.  Then $\tilde{a}=(1-x-y)a+4x$, $\tilde{b}=(1-x-y)b$, and $\tilde{c}=(1-x-y)c+4y$. The preference for $\alpha$ and $\eta$ is on account of the resulting parameter space $(\alpha,\eta)\in[-1,1]^2$.  This corresponds to $(x,y)\in\mathbb{R}^2\vert x,y\leq1\text{ and }x+y\leq1$.}
 are given in terms of the dimensionless parameters by

\begin{align}
	\tilde{a}&=\frac{4A}{E}=(1+\tilde{\delta})(1+\tilde{\sigma})=\frac{1-\alpha}{1+\alpha}(1+\delta)(1+\sigma)+4\frac{\alpha+\eta}{1+\alpha} \\
	\tilde{c}&=\frac{4C}{E}=(1+\tilde{\delta})(1-\tilde{\sigma})=\frac{1-\alpha}{1+\alpha}(1+\delta)(1-\sigma)+4\frac{\alpha-\eta}{1+\alpha} \\
	a&=\frac{4\bar{A}}{\bar{E}}=(1+\delta)(1+\sigma)\, ,\quad	c=\frac{4\bar{C}}{\bar{E}}=(1+\delta)(1-\sigma) \\
	\tilde{b}&=\frac{4B}{E}=1-\tilde{\delta}=\frac{1-\alpha}{1+\alpha}(1-\delta)\, ,\quad b=\frac{4\bar{B}}{\bar{E}}=1-\delta \,.
\end{align}

In keeping with our terminology of kinetic, centrifugal and gravitational coefficients, all of $(\tilde{\delta},\tilde{\sigma},\alpha,\eta)$ are called kinetic parameters, $(\delta,\sigma)$ are centrifugal parameters, while $\chi$ is the gravitational parameter.  The parameter $Q$ indicates the relative strength of different parts of the potential; it is related to the ratio of centrifugal to gravitational forces, and $Q=0$ gives the standard (non-rotating) double pendulum.

\subsection{RDP Parameter Space}\label{pspace}

The Lagrangian (\ref{rdp.lag}) has 8 dimensional coefficients, all energies, but they are constrained by the condition $E=\bar{E}+K$, so that in fact parameter space is 7-dimensional.  Then dimensional analysis says we should employ one dimensional parameter with units of energy (we take $E$), and 6 additional dimensionless parameters.  We take $Q$ as one of these.  We have the 5 remaining parameters in various combinations of (\ref{delta}-\ref{chi}).

\begin{figure}
	\includegraphics[scale=0.4]{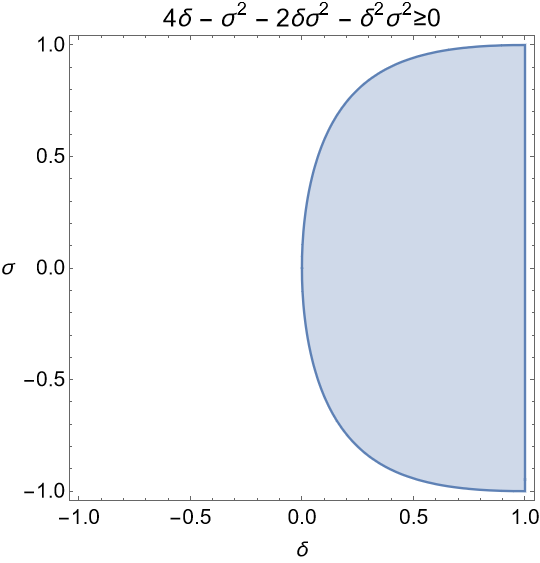}
	\caption{The `positive definite' region of the parameter space $(\delta,\sigma)\in[-1,1]^2$.}
	\label{fig.pos-def}
\end{figure}

The coefficients (\ref{coeffs}) inform the values we allow for the dimensional parameters.  Most are straightforward: $M_1,M_2,L,L_2,g,I^{(1),(2)}_{N,\perp}$ are either positive or non-negative, and we can say non-negative by including limiting cases. However, requiring that the coefficients of the Lagrangian are non-negative gives the following conditions: $\bar{A},\bar{C}\geq0$ give $0\leq I^{(1)}_P\leq M_1L_1^2+M_2L^2+I^{(1)}_\perp$ and $0\leq I^{(2)}_P\leq M_2 L_2^2+I^{(2)}_\perp$, while $K_1\geq0$ allows $L_1\geq-M_2L/M_1$.  Additionally, observe that with masses and moments non-negative, $A$ and $C$ are manifestly non-negative, so that the kinetic coefficients $A,B,C$, which are best understood as the entries of a symmetric $2\times2$ matrix, defining a quadratic form.  This quadratic form is positive definite since $AC>B^2$ (in fact, $AC-B^2\geq M_1M_2L_1^2L_2^2$).  In general, this does not follow for the centrifugal coefficients $\bar{A},\bar{B},\bar{C}$.  In terms of the dimensionless parameters, $(\tilde{\delta},\tilde{\sigma})$ are constrained by the following condition

\begin{equation}
	\tilde{a}\tilde{c}-\tilde{b}^2\,=\,4\tilde{\delta}-\tilde{\sigma}^2(1+\tilde{\delta})^2\geq0\,.\label{pos.def}
\end{equation}

\noindent but $(\delta,\sigma)$ are not likewise constrained.  This region is shown in Fig~\ref{fig.pos-def}.  If we use parameters $(\delta,\sigma,\alpha,\eta)$, then for any given values $(\delta,\sigma)\in[-1,1]^2$,  (\ref{pos.def}) is a corresponding condition on the values of $(\alpha,\eta)$.

\subsection{Special Case of Point Masses on Massless Rods (PMMR)}\label{sPMMR}

We also consider a special case of particular interest: the textbook model of point-masses $M_1$, $M_2$ on massless rods. In this case the pivot $\mathcal{O}_2$ coincides with $\mathit{CM}_{\!1}$, the masses constrained to remain distances $L_1$, $L_2$ apart from the pivots $\mathcal{O}_1$, $\mathcal{O}_2$, respectively.  This case is defined by the following relations:  $L=L_1$, and all principle moments $I^{(1),(2)}_{P,\perp,N}=0$.  However, we call this case \textit{Strict PMMR}.  There are wider classes of models that are not point-masses on massless rods, but they nevertheless exhibit identical dynamics for the subset of the parameter space where they overlap.  Observe that the Strict PMMR case has $\alpha,\eta=0$, so that the pairs $(\delta,\sigma)$ and $(\tilde{\delta},\tilde{\sigma})$ coincide.  Both pairs are thus constrained by positive definiteness, as in equation (\ref{pos.def}).  Meanwhile, the condition $L=L_1$ establishes the following common ratio between the centrifugal and gravitational coefficients

\begin{equation}
	\frac{\bar{A}}{\bar{B}}=\frac{K_1}{K_2},\quad\text{equivalently}\quad \frac{a}{b}=\frac{1+\chi}{1-\chi}\label{PMMR}
\end{equation}

We take the resulting conditions in the dimensionless parameters:  $\alpha=0$, $\eta=0$ and equation (\ref{PMMR}) to define the \textit{Strong PMMR} case.  The vanishing kinetic parameters $\alpha=\eta=0$ imply $A\bar{B}=\bar{A}B$ and $C\bar{B}=\bar{C}B$, which simplify to $A\omega_a^2=\bar{A}\omega^2$ and $C\omega_a^2=\bar{C}\omega^2$, which in turn require $I^{(i)}_\perp=I^{(i)}_N+I^{(i)}_P$ for $i=1,2$.  Thus any RDP consisting of physical pendula with principle moments satisfying these conditions will exhibit identical dynamics to the model of point masses on massless rods, provided that (\ref{PMMR}) also holds.

Furthermore, we take only the condition (\ref{PMMR}) to define the \textit{Weak PMMR} case.  It is a codimension-1 hypersurface of the (dimensionless) parameter space, or 5-dimensional, whereas Strict and Strong PMMR are 3-dimensional.  However, these three cases are nearly equivalent for questions of equilibrium and bifurcation, which do not involve the kinetic parameters (see \S\ref{eqns}, lines (\ref{equib},\ref{bif2})).  The only difference is that in Strict/Strong PMMR ($\alpha=\eta=0$), the centrifugal parameters $(\delta,\sigma)$ are still restricted to be positive definite as in (\ref{pos.def}), while in Weak PMMR this condition is relaxed.

The Weak PMMR condition (\ref{PMMR}) allows us to express the gravitational parameter $\chi$ as a rational function of $(\delta,\sigma)$ 

\begin{equation}
	\chi_{_\textit{PMMR}}=\frac{2\delta+\sigma+\delta\sigma}{2+\sigma+\delta\sigma}\,.\label{chiPMMR}
\end{equation}

\section{Equations of Motion, Equilibrium, Normal Mode Frequencies, and Bifurcation}\label{eqns}

The reduced Lagrangian has the potential function

\begin{align}
	V(\theta,\varphi)=-\frac{Q}{2}&\left(a\sin^2\theta+2b\sin\theta\sin\varphi+c\sin^2{\varphi}\right) \\
	-2&\,(1-Q)\left[(1+\chi)\cos\theta+(1-\chi)\cos\varphi\right] \nonumber
\end{align}

\noindent and the equations of motion (given be Euler-Lagrange equations) simplify to 

\begin{equation}
	\underbrace{\left[\begin{array}{cc}
			\tilde{a} & \tilde{b}\cos(\varphi-\theta) \\ \tilde{b}\cos(\varphi-\theta) & \tilde{c}
		\end{array}\right]}_{\tilde{\mathcal{M}}(\varphi-\theta)}
	\left[\begin{array}{c}
		\ddot{\theta} \\ \ddot{\varphi}
	\end{array}\right]+\tilde{b}\sin(\varphi-\theta)\left(\begin{array}{r}
		-\dot{\varphi}^2 \\ \dot{\theta}^2
	\end{array}\right)+\nabla V(\theta,\varphi)=0 \label{eom}
\end{equation}

\noindent so that equilibria of the system are extrema of $V$.  Using simplified $q=Q/(1-Q)$, (which is $q\geq0$ for $Q\in[0,1]$),

\begin{equation}
	\nabla V(\theta,\varphi)=0\quad\iff\quad\begin{array}{l}
		q(a \sin\theta+b\sin\varphi)\cos\theta-2(1+\chi)\sin\theta=0 \\
		q(b\sin\theta+c\sin\varphi)\cos\varphi-2(1-\chi)\sin\varphi=0\,.
	\end{array}\label{equib}
\end{equation}

Next, we expand the equations of motion near an equilibrium point $(\theta_0,\varphi_0)$ satisfying (\ref{equib}). Substituting $\left(\theta(\tau),\varphi(\tau)\right)=\left(\theta_0,\varphi_0\right)+\epsilon\bf{y}$ for some $\epsilon\ll1$ and $\bf{y}\in\mathbb{R}^2$, to leading order in $\epsilon$ the equations of motion (\ref{eom}) become

\begin{equation}
	\tilde{\mathcal{M}}(\varphi_0-\theta_0)\,\ddot{\bf{y}}+H[V](\theta_0,\varphi_0)\,\bf{y}=0 \label{linear}
\end{equation}

\noindent where $H[V]$ is the second derivative matrix of the potential.  From (\ref{linear}), we can examine the normal mode behavior near equilibria in the linear approximation.  Assume a common exponential behaviour

\begin{equation}
	\bf{y}(\tau)=\bf{v}e^{\pm i\sqrt{\Omega}\tau}
\end{equation} 

\noindent where $\bf{v}\in\mathbb{R}^2$ is constant, $\Omega>0$ indicates oscillatory solutions while $\Omega<0$ gives rise to two exponential modes---one growing, the other decaying.\footnote{Since both matrices in (\ref{nmf}) are symmetric, and $\tilde{\mathcal{M}}$ is positive definite (see \S\ref{pspace}), we know that (\ref{nmf}) only has real solutions.}  In (\ref{linear}) this gives 

\begin{equation}
	\left\lbrace\tilde{\mathcal{M}}(\varphi_0-\theta_0)\,(-\Omega)+H[V](\theta_0,\varphi_0)\right\rbrace\bf{v}=0\,.
\end{equation}

\noindent Thus we have the following generalized eigenvalue equation for (the square of) normal mode frequencies $\Omega$

\begin{equation}
	\det\left[-\Omega\,\tilde{\mathcal{M}}(\varphi_0-\theta_0)+H[V](\theta_0,\varphi_0)\right]=0\,.\label{nmf}
\end{equation}

\noindent  This is a quadratic in $\Omega$.  Bifurcation is precisely the scenario in which at least one root of (\ref{nmf}) is $\Omega=0$, which occurs if and only if the constant term of (\ref{nmf}) vanishes.  Thus bifurcation is given by the conditions

\begin{gather}
	\nabla V(\theta_0,\varphi_0)=0 \label{bif1} \\ 
	\det\Big( H[V](\theta_0,\varphi_0)\Big)=0 \label{bif2}
\end{gather}

\noindent These equations involve the coefficients $a,b,c$ (or equivalently, $\delta$ and $\sigma$) and parameters $q$, $\chi$.  Of course, they also involve the trigonometric functions $\cos\theta_0$, $\sin\theta_0$, $\cos\varphi_0$, $\sin\varphi_0$.  However, if these are replaced with polynomial variables, respectively $c_1$, $s_1$, $c_2$, $s_2$, and Pythagorean identities $c_1^2+s_1^2-1=0$, $c_2^2+s_2^2-1=0$ are added to the equations (\ref{bif1}-\ref{bif2}), the result is a system of polynomial equations that describe bifurcation.

The system is 5 equations in 4 variables and 4 parameters, fully given in expanded form in Appendix~\ref{outputs1}.  In Section~\ref{elim} we employ various methods of elimination for systems of polynomial equations, to eliminate the trig variables $c_1$, $s_1$, $c_2$, $s_2$ and produce a single condition on the parameters, representing the codimension-1 subset of parameter space for which the system is in bifurcation.  This is a large polynomial system, and a challenging computation for many algorithms in different implementations.  \texttt{DixonEDF} \cite{lewis2017}, written in the \textit{Fermat} computer algebra system \cite{fermat}, is a very powerful elimination algorithm, extracting factors of the resultant while computing the determinant of the Dixon matrix.

\section{Trivial Equilibria and Expected \texttt{DixonEDF} Factors}\label{trivS}

Observe that the equations for equilibrium (\ref{equib}) have constant solutions $s_1,s_2=0$, for all parameter values.  The Pythagorean identities give $c_1,c_2=\pm1$.  These configurations are $\theta,\varphi=0\text{ or }\!\pm\!\pi$; 4 combinations in total.  We refer to them as:  down-down, $\theta=\varphi=0$; down-up, $\theta=0,\varphi=\pi$; up-down, $\theta=\pi,\varphi=0$; up-up, $\theta=\varphi=\pi$, and collectively as the \textit{trivial equilibria}.  They are the only equilibria of the standard (non-rotating) double pendulum, and they are equilibria of the RDP for all parameter values.  This means questions of bifurcation, or even of the normal mode frequencies, are simple to evaluate for these equilibria, substituting the coordinates into  (\ref{bif2}) and (\ref{nmf}), respectively.

When utilizing polynomial elimination algorithms in applications, the resulting polynomial (whether a resultant, the first polynomial of a Groebner basis, etc.) often has many factors, some of them spurious.  Most algorithms or implementations compute this polynomial to full expansion, so the effective multiplying-out of all these factors adds to the computational costs.  Typically, one factor in particular, or perhaps a handful of factors, are relevant to the problem, and so the resultant polynomial wants to be factored anyway, which may itself be costly for sufficiently large outputs.  One of the comparative advantages of \texttt{DixonEDF}, then, is that as the determinant of the Dixon matrix is being computed, common numerators and denominators of rows and columns are extracted as the computation proceeds.  The end result is a list of factors of the resultant.  The factorization problem is not completely resolved, as the factors in the list are not necessarily irreducible, but the problem is often greatly reduced.  The relevant factor or factors are typically easy to identify as the longest (with the highest number of terms).

In the case of the rotating double pendulum, the trivial equilibria actually provide us with factors to expect from elimination computations -- namely, $\theta_0,\varphi_0=0\text{ or }\pi$ substituted into (\ref{bif2}) or (\ref{nmf}) for bifurcation or normal mode frequencies, respectively.  In \S\ref{elim}, we confirm precisely these factors in elimination computations.

\subsection{Trivial Bifurcations in the PMMR Case}

\begin{figure}[!ht]
	\centering
	\includegraphics[scale=.16]{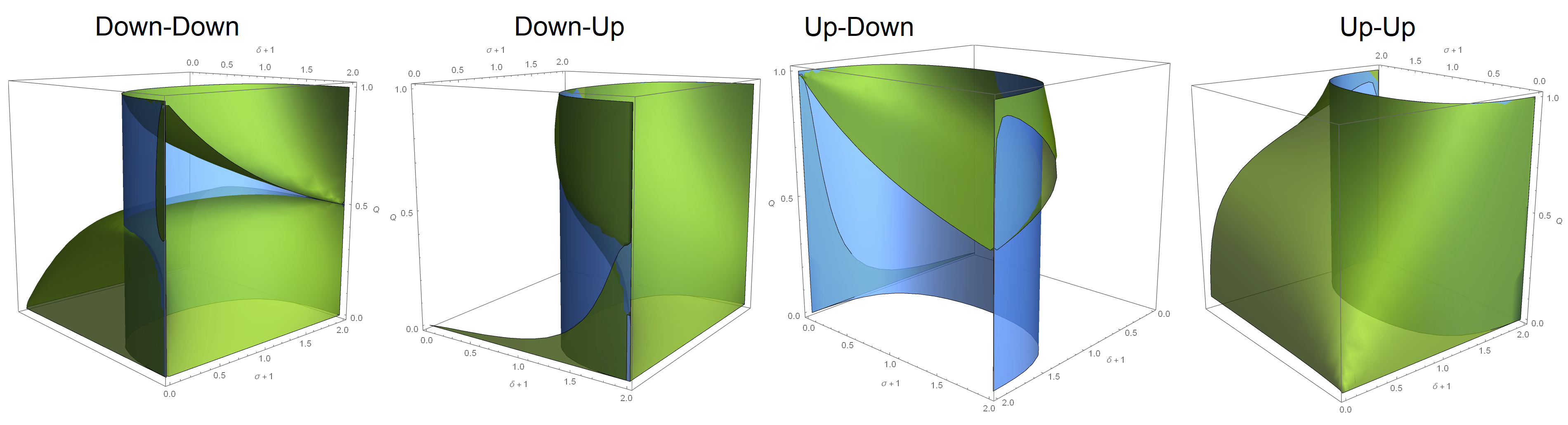}
	\caption{Contour Plots of bifurcation condition: $\det\left(H[V](\theta,\varphi)\right)=0$, evaluated at the four trivial equilibria.  The horizontal axes measure $1+\delta$ and $1+\sigma$, while the vertical is $Q$.  The vertical blue surface is the same in each plot---it is the surface where the centrifugal parameters $(\delta,\sigma)$ form a singular quadratic form, ie. equality in (\ref{pos.def}).  Thus the positive-definite region is the convex one, with only positive values of $\delta$, corresponding to Strong PMMR, whereas Weak PMMR allows the entire box $[-1,1]^2\times[0,1]$.  Thus we may observe that the up-up equilibrium never undergoes bifurcation under Strong (or Strict) PMMR, while it does in Weak PMMR.}
	\label{trivbif}
\end{figure}

In any of the PMMR special cases introduced in \S\ref{sPMMR}, elimination of the polynomial system for bifurcation, with (\ref{chiPMMR}) substituted in for $\chi$ and numerators taken, will result in a single polynomial condition on the three parameters $(\delta,\sigma,Q)$.  We thus anticipate visualizing the degenerate models in parameter space with 3D contour plots.  Likewise, the expected factors given by the trivial equilibria are polynomial conditions in the same space, indicating for which RDP models is, say for example, the down-down equilibrium in bifurcation.  Figure~\ref{trivbif} shows the bifurcations of the four trivial equilibria in the PMMR subcase.

The translucent blue surfaces in each plot of Fig.~\ref{trivbif} are the same---the surface given by equality of (\ref{pos.def}) for $(\delta,\sigma)$.  Thus the convex region on one side of the blue surface, with only positive values of $\delta$, corresponds to positive definite centrifugal coefficients, matching the Strong PMMR case.  The entire space $(\delta,\sigma,Q)\in[-1,1]^2\times[0,1]$ is allowed by Weak PMMR.

The bifurcation equation $~\det\left(H[V](\theta,\varphi)\right)=0~$ is quadratic in $Q$, so formally we can write two solutions as the separate branches of the quadratic formula.  We observe that for the down-down equilibrium, both branches give physical values $Q\in[0,1]$, although the ``first" of these, with smaller $Q$-value, occurs in Weak PMMR, for any values $(\delta,\sigma)\in[-1,1]^2$, whereas the second branch only exists in Strong PMMR.  For the down-up and up-down equilibria, the lower branch gives $Q\in[0,1]$ across the positive-definite region, but both branches have solutions for only very narrow regions outside of positive-definite.  Perhaps most strikingly, only one of the branches for the up-up equilibrium presents a bifurcation, and that only outside the positive-definite region.

\section{Elimination Computations for the Rotating Double Pendulum}\label{elim}

In this section we report the results of various elimination computations on polynomial systems derived from the rotating double pendulum, including performance comparison of \texttt{DixonEDF} to other techniques in the software systems \texttt{Maple} and \texttt{Magma}.  

\subsection{PMMR Equilibria}

First, to get a sense for what equilibria exist for the RDP model, and what bifurcations to expect, we compute resultants of the equilibrium system.  It should be noted that we think of bifurcation primarily in terms of varying $Q$; the double pendulum rotating faster or slower.  This also corresponds to the fact that the dimensionless parameters, apart from $Q$, depend on how the pendulum is constructed, whereas $Q$ is the only parameter that characterizes the environment which the RDP is placed in: namely, the ratio of the strength of gravity to the strength of the centrifugal forces due to rotation.

Early in this study, we only considered the Strict PMMR case (generalization to physical pendula came later), in which the potential function can be parameterized in terms of 3 parameters; in addition to $Q$, the ratio of the pendulum masses $M_2/M_1$ and ratio of lengths $L_2/L_1$.  To sample these dimensions of parameter space, we chose 3 values for each:  $M_2/M_1 = 2/3,\, 1$, or $3/2$,  and $L_2/L_1=3/4,\, 1$, or $4/3$.  The combinations of these give the values for $(\delta,\sigma,\chi_\text{PMMR})$ in Table~\ref{OG9tab}.

\begin{table}
	\centering
	\begin{tabular}{|l|l|l|l|}
		\hline
		& $L_2/L_1=3/4$ &  $L_2/L_1=1$  & $L_2/L_1=4/3$  \\
		\hline
		$M_2/M_1=2/3$ & $\delta=\tfrac{25}{73},~\sigma=\tfrac{31}{49},~\chi=\tfrac{7}{13}\vphantom{\Big\vert}$
		& $\delta=\tfrac{3}{11},~\sigma=\tfrac{3}{7},~\chi=\tfrac{3}{7}$
		& $\delta=\tfrac{29}{125},~\sigma=\tfrac{13}{77},~\chi=\tfrac{7}{23}$ \\
		\hline
		$M_2/M_1=1$   & $\delta=\tfrac{17}{65},~\sigma=\tfrac{23}{41},~\chi=\tfrac{5}{11}\vphantom{\Big\vert}$
		& $\delta=\tfrac{1}{5},~\sigma=\tfrac{1}{3},~\chi=\tfrac{1}{3}$
		& $\delta=\tfrac{5}{29},~\sigma=\tfrac{1}{17},~\chi=\tfrac{1}{5}$ \\
		\hline
		$M_2/M_1=3/2$ & $\delta=\tfrac{35}{179},~\sigma=\tfrac{53}{107},~\chi=\tfrac{11}{29}\vphantom{\Big\vert}$
		& $\delta=\tfrac{1}{7},~\sigma=\tfrac{1}{4},~\chi=\tfrac{1}{4}$
		& $\delta=\tfrac{7}{55},~\sigma=-\tfrac{1}{31},~\chi=\tfrac{1}{9}$ \\
		\hline
	\end{tabular}
	\caption{Values of $\delta,~\sigma$ and $\chi$ for 3 values each of the mass ratio and length ratio of a rotating double pendulum consisting of point-masses on massless rods.}\label{OG9tab}
\end{table}

To visualize the non-trivial equilibria, depending on $Q$, for RDP parameters in Table~\ref{OG9tab}, we will use an alternative to the polynomial system made by inclusion of Pythagorean theorems.  Rather, we will parameterize the angles $\theta$ and $\varphi$ by their half-tangents.

\begin{align}
	t=\tan\left(\theta/2\right)\quad\implies\quad&\sin\theta=\frac{2t}{1+t^2},\quad\cos\theta=\frac{1-t^2}{1+t^2} \label{hatp1}\\
	s=\tan\left(\varphi/2\right)\quad\implies\quad&\sin\varphi=\frac{2s}{1+s^2},\quad\cos\varphi=\frac{1-s^2}{1+s^2} \label{hatp2}
\end{align}

With these rational parameterizations, the equations for equilibrium  $\nabla V(\theta,\varphi)=0$ become

\begin{align}
	\frac{\partial V}{\partial\theta}=0\quad\implies\quad &q s - 2 t + q t - 2 s^2 t + q s^2 t - 2 t^3 - q t^3 - 2 s^2 t^3 - q s^2 t^3 - q s t^4 - q s \delta + q t \delta \nonumber \\ 
	&+ q s^2 t \delta - q t^3 \delta - q s^2 t^3 \delta + q s t^4 \delta + q t \sigma + q s^2 t \sigma - q t^3 \sigma - q s^2 t^3 \sigma + q t \delta \sigma + q s^2 t \delta \sigma \nonumber \\
	& - q t^3 \delta \sigma - q s^2 t^3 \delta \sigma - 2 t \chi - 2 s^2 t \chi - 2 t^3 \chi - 2 s^2 t^3 \chi =0 \label{poly1}\\
	\frac{\partial V}{\partial\varphi}=0\quad\implies\quad &-2 s + q s - 2 s^3 - q s^3 + q t - q s^4 t - 2 s t^2 + q s t^2 - 2 s^3 t^2 - q s^3 t^2 + q s \delta - q s^3 \delta \nonumber \\
	&- q t \delta + q s^4 t \delta + q s t^2 \delta - q s^3 t^2 \delta - q s \sigma + q s^3 \sigma - q s t^2 \sigma + q s^3 t^2 \sigma - q s \delta \sigma+ q s^3 \delta \sigma \nonumber \\
	& - q s t^2 \delta \sigma + q s^3 t^2 \delta \sigma + 2 s \chi + 2 s^3 \chi + 2 s t^2 \chi + 2 s^3 t^2 \chi =0 \label{poly2}\\
	&\text{These polynomials are the numerators of the partial derivates, when sines } \nonumber \\
	&\text{and cosines are replaced by the half-angle tangent parametrizations (\ref{hatp1},\ref{hatp2}).} \nonumber
\end{align}

We compute three resultants of the polynomials (\ref{poly1}) and (\ref{poly2}), eliminating in turn each of $t,~s$ and $q$.  This is done very simply with the built-in Mathematica command \texttt{Resultant}, which may be used to eliminate one variable from a system of two polynomial equations.

\begin{align}
	f(s,q,\delta,\sigma,\chi)&=\text{Largest irreducible factor  of }
	~\texttt{Resultant}\left[\left\lbrace(\ref{poly1}),(\ref{poly2})\right\rbrace,t\right]\label{resF} \\
	g(t,q,\delta,\sigma,\chi)&=\text{Largest irreducible factor  of }
	~\texttt{Resultant}\left[\left\lbrace(\ref{poly1}),(\ref{poly2})\right\rbrace,s\right]\label{resG} \\
	h(t,s,\delta,\sigma,\chi)&=\text{Largest irreducible factor  of }
	~\texttt{Resultant}\left[\left\lbrace(\ref{poly1}),(\ref{poly2})\right\rbrace,q\right]\label{resH}
\end{align}

\begin{figure}
	\centering
	\includegraphics[scale=0.28]{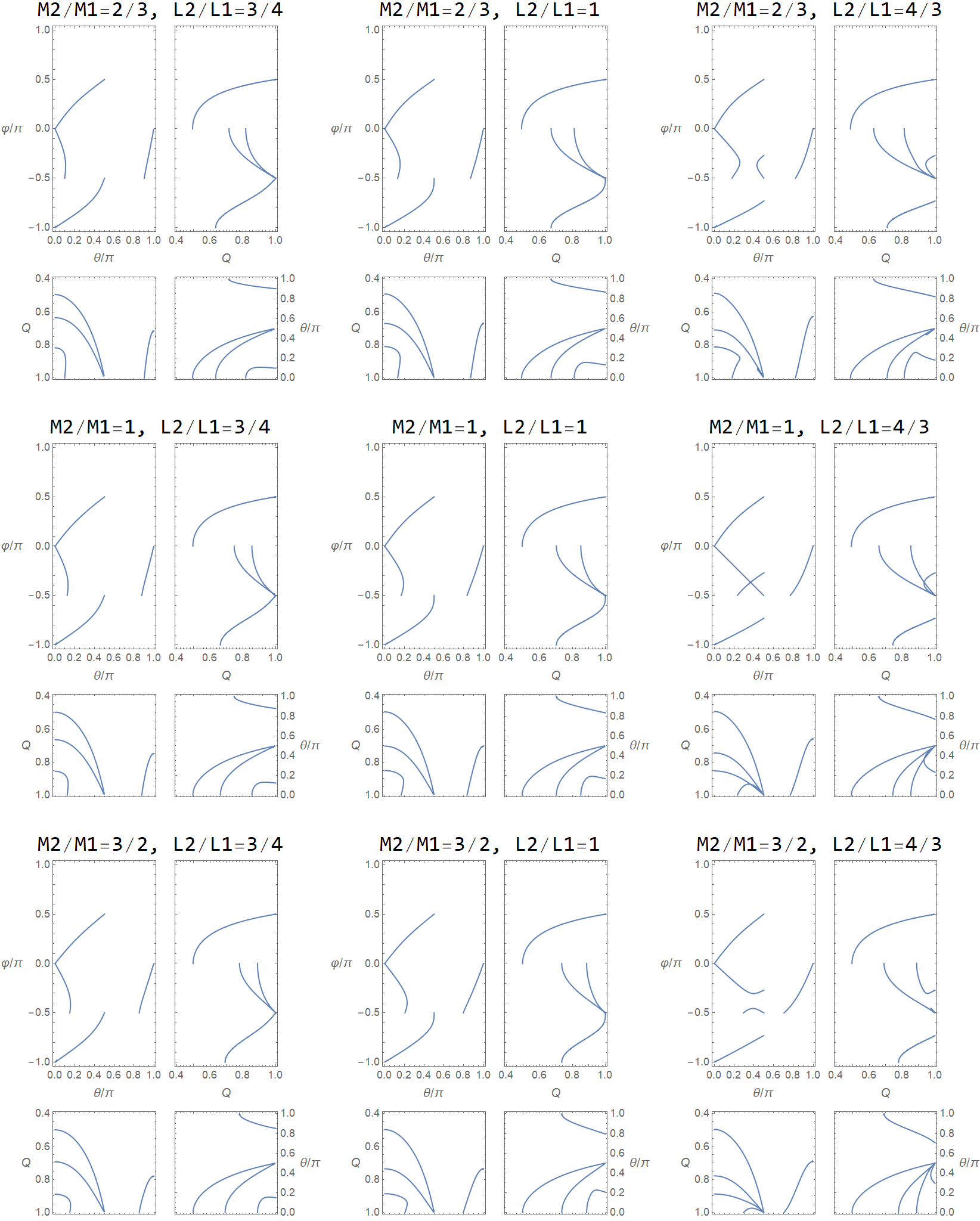}
	\caption{For each of the mass and length ratio pairs in Table~\ref{OG9tab}, this figure visualizes the equilibria of the RDP---solutions of (\ref{equib})---across continuous values of the parameter $Q$, demonstrating both trivial and non-trivial equilibria.}
	\label{OG9proj}
\end{figure}

The largest irreducible factors of the resultants are of interest; call these $f,~g$ and $h$.  The other factors are products of powers of $t,~s,~1+t^2$ and $1+s^2$.  The factors $f$ and $h$ are both 1290 terms, while $h$ is 32 terms. Of course, $h$ is very easy to arrive at; simply solving both (\ref{poly1}),(\ref{poly2}) for $q$---they are linear in $q$---setting the results equal, collecting everything to one side on common denominator, and take the numerator.  For the parameter values in Table~\ref{OG9tab}, Fig.~\ref{OG9proj} shows the 0-contours of these resultants, giving projections of the equilibria into $(\theta,\varphi)$, $(Q,\theta)$ and $(Q,\varphi)$ planes for $\theta\geq0$ and $Q\in[0,1]$.

For the 9 sets of parameter values, each subplot of Fig.~\ref{OG9proj} consists of 4 panels.  The top-left panel shows the equilibria---solutions of (\ref{equib})---in $(\theta,\varphi)$-space, for all $Q$-values $\in[0,1]$.  Curves coming from either the corners of the plots or the $\theta$-axis at $\theta=0,~\pi$ indicate equilibria originating at the trivial bifurcations (see \S\ref{trivS}), while curves elsewhere indicate additional, non-trivial bifurcations.  The top-right panel shows the solutions $\varphi(Q)$ of (\ref{equib}) \textbf{for $\theta>0$ only}.  The bottom panels are the same plot, but rotated 90$^\circ$ relative to each other, and shows the solutions $\theta(Q)>0$ of (\ref{equib}).  In the bottom-left the $Q$-axis descends down the vertical, the $\theta$-axis coinciding with that of the top-left panel above.  In the bottom-right, the $Q$-axis is the horizontal, matching that of the $\varphi(Q)$-plot at top-right.  The top-left panels are produced by a transformation of 0-contours $h=0$ of the resultant (\ref{resH}), while the top-right panels are transformations of 0-contours $f=0$ of the resultant (\ref{resF}).  The bottom panels are similarly produced by a transformation of $g=0$ of the resultant (\ref{resG}).

In combination, these panels allow to be read off which $\varphi(Q)$ and $\theta(Q)$ solutions pair together as equilibria in the top-left, as well as the bifurcation $Q$-values at which solutions diverge.  In particular, we observe non-trivial bifurcations, typically associated with two solutions of (\ref{equib}), for sufficiently large $Q$, in the vicinity of $(\theta=\pi/2,~\varphi=-\pi/2)$.  We thus anticipate at least one additional relevant factor in the results of a potential \texttt{DixonEDF} computation for the bifurcation system (\ref{bif1},\ref{bif2}), beyond the four expected factors corresponding to the trivial bifurcations.   Within this 9-point sample of the remaining parameter space, we also notice that the non-trivial bifurcations only occur for $L_2/L_1=4/3$.  

\subsection{\texttt{DixonEDF} Computation of the RDP Bifurcation Resultant}

Bifurcation of the RDP system is described by the system of equations (\ref{bif1},\ref{bif2}).  Taken together with Pythagorean identities to constrain trigonometric variables, this gives a system of 5 polynomial equations.  When expressed in terms of:  $c_1,~s_1,~c_2,~s_2,~q=\nicefrac{Q}{1-Q},~d=1+\delta,~s=1+\sigma$ and $\chi$, these equations have 5, 6, 48, 3 and 3 terms, respectively.  These are given in full in Appendix~\ref{outputs1}.  The 48 terms, corresponding to the bifurcation equation (\ref{bif2}), is at most cubic in trigonometric variables, quadratic in parameters, and has maximum total order 10.

The \textit{Fermat} implementation of \texttt{DixonEDF} succeeds in computing the resultant for this system, eliminating $c_1,~s_1,~c_2$ and $s_2$.  Specifically, \texttt{DixonEDF} extracts common factors while computing the determinant of the Dixon matrix, so that the computation identifies factors with the following number of terms.

\begin{equation}
	\text{Lengths of factors of the Dixon resultant:}~ 1 ~1 ~5 ~2 ~1 ~23288 ~2 ~2 ~1 ~1
\end{equation}

\noindent However, the 23288-term polynomial does factor.  Four expected factors, corresponding to bifurcations of the constant equilibria $\theta,\varphi=0\text{ or }\!\pm\pi$, are found by evaluating the 48-term polynomial with $s_1,s_2=0,~c_1,c_2=\pm1$, which gives the following 10-term polynomials

\begin{align}
	\text{down-down: }&c_1=c_2=1 \nonumber \\
	4 - 4 \chi^2& - 4 d q - 4 \chi  d q - 4 q^2 + 4 d q^2 - d^2 q^2 + 4 \chi d q s + 2 d^2 q^2 s - d^2 q^2 s^2 \label{dd}\\
	\text{down-up: }&c_1=1,~c_2=-1 \nonumber \\
	-4 + 4 \chi^2& - 4 d q - 4 \chi d q - 4 q^2 + 4 d q^2 - d^2 q^2 + 4 d q s + 2 d^2 q^2 s - d^2 q^2 s^2 \\
	\text{up-down: }&c_1=-1,~c_2=1 \nonumber \\
	-4 + 4 \chi^2& + 4 d q + 4 \chi d q - 4 q^2 + 4 d q^2 - d^2 q^2 - 4 d q s + 2 d^2 q^2 s - d^2 q^2 s^2 \\
	\text{up-up: }&c_1=c_2=-1 \nonumber \\
	4 - 4 \chi^2& + 4 d q + 4 \chi d q - 4 q^2 + 4 d q^2 - d^2 q^2 - 4 \chi d q s + 2 d^2 q^2 s - d^2 q^2 s^2 \label{uu}
\end{align}

The 23288-term factor is divisible precisely by the polynomials (\ref{dd}-\ref{uu}).  Initializing the list of denominators with these expected factors, the following list is then found

\begin{equation}
	\text{Lengths of factors, taking account of (\ref{dd}-\ref{uu}):}~ 5 ~1 ~2 ~2 ~1 ~2 ~1 ~2 ~2 ~2 ~5 ~5 ~1 ~1 ~5 ~3 ~3 ~175~ 6744\label{6744}
\end{equation}

\noindent where the 175-term factor is the product of (\ref{dd}-\ref{uu}).  Thus we are most interested in the 6744-term polynomial, which contains all of $(q,d,s,\chi)$.  Many of the smaller factors depend on a strict subset of these parameters.  Repeating the computation using alternate parameters to $d$ and $s$, the longest factor is then 3257 terms.  Performing the change-of-variables to our 6744-term factor matches the 3257 terms, a confirmation that these results are at least consistent.  The computation resulting in (\ref{6744}) took 16 seconds and 134 MB of RAM on a 24-core Intel Xeon Gold 6126 CPU running Linux Redhat 6 with access to 150 GB of RAM and 2.3 TB of storage.  The computation with alternate parameters takes only 7.24 seconds, although 139 MB.  Attempts with the \texttt{FGb} \textit{Maple} package for Grobner bases, and a Groebner basis method in \textit{Magma} by Allan Steel\footnote{Both of these Groebner basis implimentations are referenced in \cite{lewis2017}} were both terminated after tens of hours and several gigabytes of RAM used, or more.

\subsection{Bifurcation in PMMR Special Case}

Considering bifurcation in the special case of point-masses on massless rods, we can proceed in two ways. We can take the 6744-term resultant from the computation described above and substitute for $\chi$ according to (\ref{chiPMMR}).  Or we can make the substitution into the bifurcation system of equations (\ref{bif1},\ref{bif2}), taking numerators (the resulting polynomials are given in Appendix~\ref{output2}), and compute a resultant.  We have completed both of these calculations, and compare the results.

\begin{figure}[h!]
	\centering
	\includegraphics[scale=0.24]{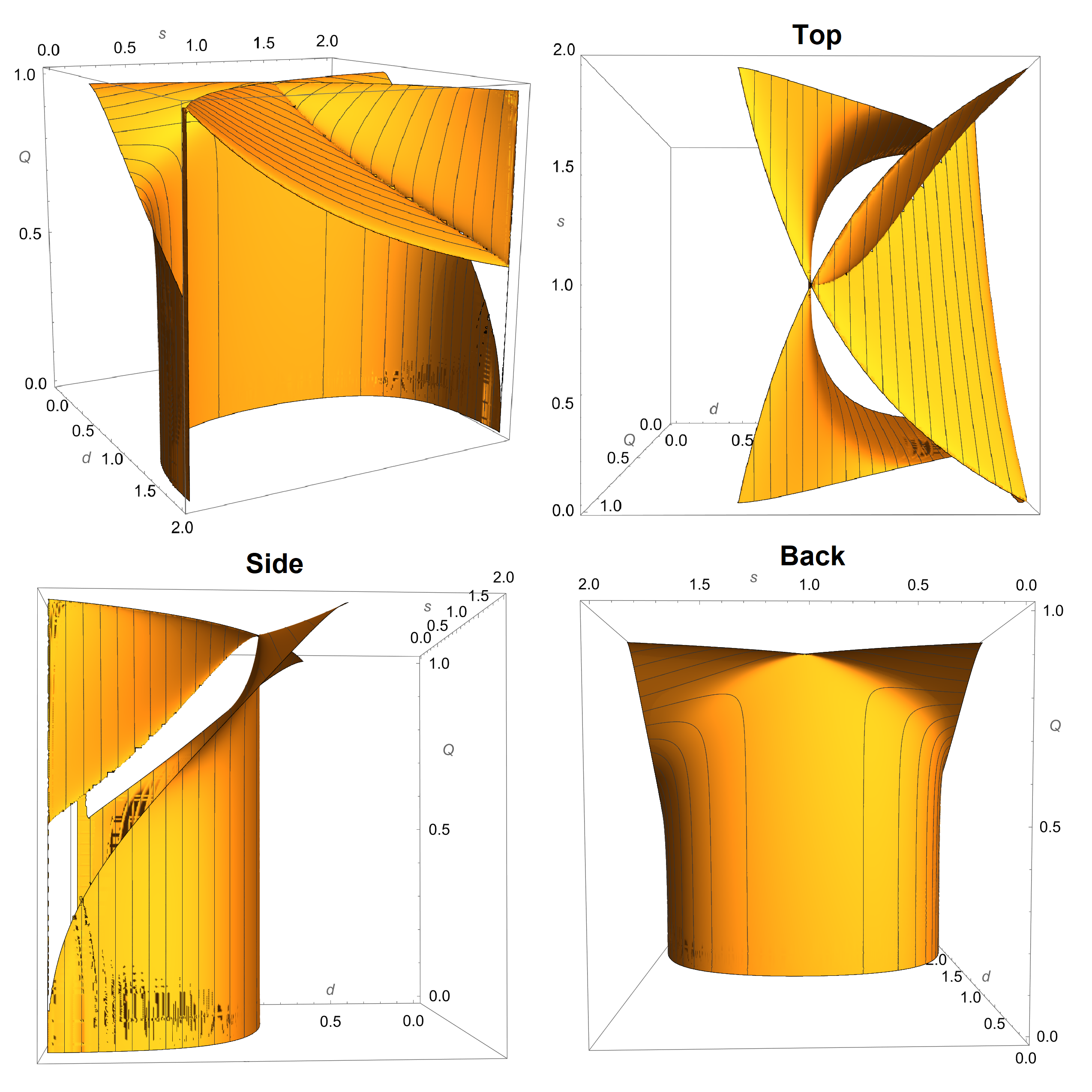}
	\caption{Contour plot of $R_{3356}=0$.  Subset of PMMR parameter space $(d,s,Q)$ which exhibits the non-trivial bifurcation identified by the 3356-term resultant.}\label{bifresCP}
\end{figure}

The \texttt{DixonEDF} computation for the PMMR bifurcation system finds factors with the following numbers of terms

\begin{equation}
	\text{Lengths of factors: }~8 ~3 ~13 ~2 ~1 ~1 ~1 ~1924 ~5 ~1 ~2 ~21 ~23 ~23 ~21,
\end{equation}

\noindent this computation using 19s and 92 MB. The 1924-term factor is precisely the result of the substitution $\chi=\nicefrac{d(1+s)-2}{2+d(s-1)}$, according to (\ref{chiPMMR}), into the longest factor of (\ref{6744}).  Replacing $q\rightarrow Q/(1-Q)$ in these 1924 terms, the primary factor of the numerator is 3356 terms long.  Fig.~\ref{bifresCP} shows a contour plot of this polynomial $R_{3356}=0$ in $(d,s,Q)$-space.

\section{Conclusion}

In this work we have found that both the equilibrium configurations and bifurcation parameter values which produce degenerate equilibria of the rotating double pendulum can be described by systems of polynomial equations.  Several questions about these features can then conceivable be answered by employing algorithms of computational algebraic geometry, in particular elimination by resultant computations. While the system of equations for bifurcation appears to be intractable for classic algorithms such as Grobner bases, as exemplified by two routines in \textit{Maple} and \textit{Magma}, the \textit{Fermat}  implementation of the \texttt{DixonEDF} \cite{fermat,lewis2017} algorithm, computing the Dixon resultant with early detection of factors, is able to complete these computations on the order of seconds.  With early detection of factors, these calculations recover the trivial bifurcations which occur at the (trivial) vertical equilibria of the RDP, as well as confirm and describe a fifth\footnote{We note here that the non-trivial bifurcation we have found for the RDP, as described by the resultant $R_{3356}=0$, may actually be rightly understood as \textit{three} additional bifurcations: there is a narrow region, approximately the `crease' of the surface in Fig. \ref{bifresCP}, within which a vertical line (particular $\delta,\sigma$ values chosen) intersects the surface three times; for three distinct $Q$-values $\in(0,1)$.  Thus three non-trival bifurcations would be expected in the corresponding bifurcations plots.  This detail must be analyzed in future work.} ~`non-trivial' bifurcation which occurs for comparatively high $Q$-values.  This work suggests that the equilibrium and bifurcation structures of other physical systems could be analyzed by the use of similar computations, provided that the equations (\ref{bif1},\ref{bif2}) can be expressed (or perhaps approximated) as systems of polynomial equations.

This work also suggests a potential link with the work in \cite{maiti2016}.  It would be conjectured that the transition to quasi-periodic behaviour that Maiti et. al. observe via numerical experiments corresponds to the first bifurcation of $(\theta=0,\varphi=0)$, which occurs near $Q\approx1/2$ for generic values of the other parameters fixed.  This bifurcation appears to be a pitchfork (this should be confirmed in future work), with the vertical configuration $(0,0)$ becoming unstable for increasing $Q$, and spawning two stable equilibria which are $\pm(\theta_*,\varphi_*)\neq0$ by the symmetry of the system.

\subsection{Further Work}

There are many directions in which this work can be continued

\subsubsection{Normal Mode Analysis}

As demonstrated in Appendix \ref{rpp}, one desires to characterize the nature of the RDP equilibria and bifurcations by completing the normal mode analysis.  This invokes the system of equations (\ref{bif1}) and (\ref{nmf})

\begin{align}
	\nabla&V(\theta_0,\varphi_0)=0 \\
	\det&\!\left[-\Omega\,\tilde{\mathcal{M}}(\varphi_0-\theta_0)+H[V](\theta_0,\varphi_0)\right]=0 \label{nmf2}
\end{align}

\noindent which reduces to the bifurcation system (\ref{bif1},\ref{bif2}) for $\Omega=0$.  Preliminary investigations already show that this problem, attempting to eliminate the angle variables, is approaching the limit of \texttt{DixonEDF} practical capabilities, at least on the hardware that has been in use.  Moreover, the resultants we have found are on the order of 1 million terms long, and thus they are impractical to use for further calculations.  However, rather than considering this system in full, observe that the second equation is a quadratic in $\Omega$.  Since both matrices in (\ref{nmf2}) are symmetric, and the kinetic matrix $\tilde{\mathcal{M}}(\varphi_0-\theta_0)$ is positive-definite, we know (\ref{nmf2}) always has real roots, and we are interested in the sign of these roots.  Thus we may consider various questions of the signs of the coefficients of this quadratic polynomial.  This line of investigation is currently underway.

\subsubsection{Other Special Cases}

\begin{figure}[h!]
	\includegraphics[scale=0.7]{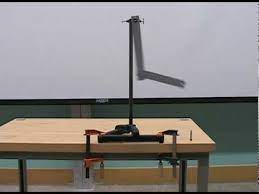}
	\caption{A double pendulum constructed from similar uniform slabs; a frame from a demonstration video by YouTube user \texttt{stevenbtroy} \cite{u.slabs}}\label{uni_slabs}
\end{figure}

A further special case that would be in some sense natural to consider, is that of the double pendulum constructed from uniform slabs, ie. rectangular prisms, made of the same material, as depicted in Fig. \ref{uni_slabs}.  This case similarly reduces the parameter space by one dimension, since the ratio $L_1/L$ is $1/2$.  At most minimal, one may consider pendula of 0 cross-section but uniform linear mass density, with the joints fixed at the end points of the first pendulum.  In this sense, the uniform slabs case is something like a counterpoint to the PMMR case.  The minimal subcase (perhaps we would refer to \textit{thin rods}) reduces the relevant parameter space a dimension further, as masses are proportional to the lengths, in addition to the kinetic parameters $\alpha,\eta$ vanishing.

\begin{figure}[h!]
	\includegraphics[scale=0.5]{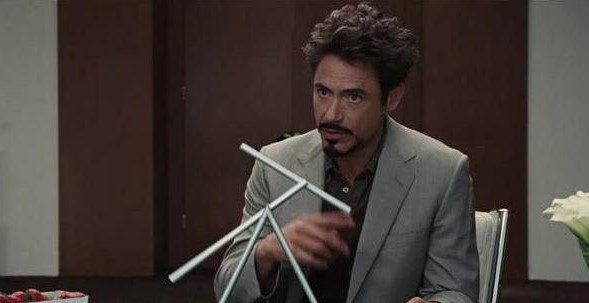}
	\caption{The kinetic sculpture Swinging Sticks${}^\text{\textregistered}$ \cite{swinging.sticks}, as seen in the film Iron Man (2008) \cite{ironman}}\label{iron_man}
\end{figure}

Yet another case which could be considered, was made somewhat famous, or perhaps infamous, by the inclusion of a certain kinetic sculpture in the original Iron Man movie (2008) \cite{ironman}, namely so-called Swinging Sticks${}^\text{\textregistered}$ \cite{swinging.sticks}, as pictured in Fig. \ref{iron_man}.  

This construction appears to include the following properties: not only is the \textit{CM} of the second pendulum very close to the second pivot, giving a comparably low $K_2$ value, furthermore the \textit{CM} of the first pendulum is behind the first pivot, corresponding to a negative value $L_1<0$ and thus a comparatively lower $K_1$ value.  If this construction does indeed keep $K_1\geq0,$\footnote{Equivalently $\chi\in[-1,1]$, and in any case if it does not, our analysis would be in this regime, at least initially.} then these properties would mean the Swinging Sticks coincide with our model of the double pendulum for larger $Q$-values when compared to the generic construction and a given rotation frequency $\omega_a$.  Thus if a swinging sticks model were employed to experimentally study the dynamics of the RDP, we may expect to be able to observe the effects of the non-trivial bifurcation, which occurs for higher $Q$-values, with moderate rotation speeds.

\subsubsection{Generalizations}

We finally mention some further generalizations that could be added to the RDP.  As many other researches have considered \cite{samaranyake1993, bridges2001}, one may take the transverse pendulum: the axes of the two joints of the pendulum not restricted to parallel.  It would be expected that rotation would tend to disturb the doubly-degenerate equilibrium that Bridges and Georgio \cite{bridges2001} study, but it would be natural to ask whether the rotation then might produce this degeneracy at different parameter values or for other equilibrium solutions.

One may also shift the pendulum horizontally from the axis of rotation, either in the plane of the pendulum or perpendicular to it.  Indeed, it should be remarked that our anaylsis here (specifically, the derivation of the Lagrangian in Appendix \ref{Ap.RDP.deriv}) assumes that all the pivots and centers of masses are co-planar with the ``plane of the pendulum", whereas physical contstructions typically have the second pendulum offset from the first in this direction.  However, it is easy to verify (for a non-transverse pendulum) that any fixed translation perpendicular to the pendulum plane, either of the first mass at the inner pivot or of the second mass at the outer pivot, has no effect on the dynamics; the terms added to the Lagrangian (\ref{rdp.lag}) by such displacements are total derivatives.  However, this generalization may be expected to have non-trivial implications if one considers a transverse double pendulum (of course in this case, a single `plane of the pendulum' doesn't exist).

If, on the other hand, the pendulum is translated horizontally in it's plane, there are non-trivial consequences.  Even just considering the rotationg physical pendulum of Appendix \ref{rpp}, shifting the rotation axis away from the inner pivot introduces much larger centrifugal forces for small rotation speeds.  The trivial equilibria are no longer equilibrium solutions for all parameter values, but must be considered functions of parameter space.  Given that the non-shifted rotating physical pendulum exhibits a pitchfork bifurcation, one might expect to see the pitchfork perturbed under this alteration of the system.

\newpage

\appendix

\section{A Single Rotating Physical Pendulum}\label{rpp}

In this Appendix, we consider a much simpler system: a single physical pendulum, made to rotate uniformly about the vertical axis through the pivot.  We find all equilibria and bifurcations of the system, and assess the normal mode frequency of the system for any non-degenerate equilibria.  This is partly to demonstrating the scope of analysis which the authors hope to complete in studying the RDP, but also a step towards the derivation of the double pendulum Lagrangian in Appendix \ref{Ap.RDP.deriv}.

Consider a single spinning physical pendulum of mass $M$, with principal axes aligned in the following way:

\begin{itemize}
	\item The pivot $\mathcal{O}$ and center of mass $\mathit{CM}$ a distance $\ell$ apart along one of the principal axes of the rigid body, with moment of inertia $I_P$ ($P$ for "pendulum", ie. the axis of the arm of the pendulum).  This axis is an angle $\theta$ from the vertical.
	\item Another principal axis, with moment $I_N$,  normal to the plane of the pendulum. This plane rotates about the vertical axis through the pivot with angular frequency $\omega_a$.
	\item The third principal axis, perpendicular to the previous two, has moment $I_{\perp}$.
	\item In a typical construction, $I_N$ might be the largest moment, and $I_P$ the smallest, but this does not necessarily have to be the case.
\end{itemize}

\begin{figure}[!ht]
	\includegraphics[scale=.5]{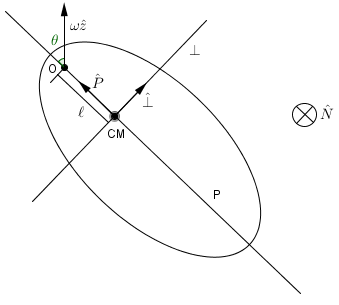}\label{rpp.fig}
	\caption{A schematic showing the orientation of a solid body as a rotating physical pendulum.  }
\end{figure}

With these axes ordered ($P,\perp,N$), set the basis vectors $\hat{P}$ pointing from $\mathit{CM}$ to the pivot, $\hat{N}$ into the page (ie. in the direction of counter-clockwise rotation when seen from above), and so $\hat{\perp}=\hat{N}\times\hat{P}$ pointing as shown in Fig. \ref{rpp.fig}, keeping a right-hand oriented frame.  In this frame the moment of inertia tensor $[I]$ is diagonal, and the angular velocity is

\begin{equation}
	\vec{\omega}=\omega_a\hat{z}-\dot{\theta}\hat{N}=\left(\omega_a\cos\theta,\omega_a\sin\theta,-\dot{\theta}\right)\,.
\end{equation}

\noindent The kinetic energy of the system is

\begin{align}
	T&=\frac{1}{2}\vec{\omega}\cdot[I]\vec{\omega}+\frac{1}{2}Mv_\mathit{CM}^2 \\
	&=\frac{1}{2}I_P\omega_a^2\left(1-\sin^2\theta\right)+\frac{1}{2}I_\perp\omega_a^2\sin^2\theta+\frac{1}{2}I_N\dot{\theta}^2+\frac{1}{2}Mv_{CM}^2 \nonumber
\end{align}

The velocity of the center of mass has two perpendicular components: one due to rotation about the vertical with angular frequency $\omega_a$, and the other due to the pendulum swinging, ie. changing $\theta$

\begin{equation}
	\vec{v}_{CM}=\ell\dot{\theta}\hat{\perp}+\left(\ell\sin\theta\right)\omega_a\hat{N}
\end{equation}

\noindent Thus $v_{CM}^2=\ell^2(\dot{\theta}^2+\omega_a^2\sin^2\theta)$, and so

\begin{align}
	T=&\cancelto{\text{ignore}}{\frac{1}{2}I_P\omega_a^2}+\frac{1}{2}I_N\dot{\theta}^2+\frac{1}{2}\left(I_\perp-I_P\right)\omega_a^2\sin^2\theta+\frac{1}{2}M\ell^2\left(\dot{\theta}^2+\omega_a^2\sin^2\theta\right) \nonumber \\
	=&\frac{1}{2}\left(I_N+M\ell^2\right)\dot{\theta}^2+\frac{1}{2}\left(I_\perp-I_P+M\ell^2\right)\omega_a^2\sin^2\theta \,.
\end{align}

The gravitational potential energy $U$, with $U=0$ at $\mathcal{O}$, is $U=-Mg\ell\cos\theta$, making the Lagrangian of the system

\begin{align}
	\mathcal{L}_1=&\frac{1}{2}\left(M\ell^2+I_N\right)\dot{\theta}^2+\frac{1}{2}\left(M\ell^2+I_\perp-I_P\right)\omega_a^2\sin^2\theta \nonumber\\
	&+Mg\ell\cos\theta\,
\end{align}
\noindent where the subscript-1 is to indicate the single pendulum; the generalization to any number of pendula linked together being $\mathcal{L}_n$ for $n$ a positive integer.

We now non-dimensionalize time, and the Lagrangian will simplify considerably. Let $t=\omega \tau$, with $\omega$ to be determined shortly.  With $\theta^\prime=d\theta/d\tau=\dot{\theta}/\omega$, the Lagrangian is

\begin{equation}
	\mathcal{L}_1=\frac{1}{2}\left(M\ell^2+I_N\right)\omega^2{\theta^\prime}^2+\frac{1}{2}\left(M\ell^2+I_\perp-I_P\right)\omega_a^2\sin^2\theta+Mg\ell\cos\theta
\end{equation}

\noindent Let $E=(M\ell^2+I_N)\omega^2$, $\bar{E}=(M\ell^2+I_\perp-I_P)\omega_a^2$, $K=Mg\ell$, and define $\omega^2$ by $E=\bar{E}+K$. That is,

\begin{equation}
	\omega^2=\frac{(M\ell^2+I_\perp-I_P)\omega_a^2+Mg\ell}{M\ell^2+I_N}\label{def.omeg}
\end{equation}

\noindent and finally let

\begin{equation}
	Q=\frac{\bar{E}}{E}=\frac{(M\ell^2+I_\perp-I_P)\omega_a^2}{(M\ell^2+I_\perp-I_P)\omega_a^2+Mg\ell}\in[0,1]\,.
\end{equation}

This gives the normalized Lagrangian

\begin{equation}
	\mathbf{L_1}=\frac{\mathcal{L}_1}{E}=\frac{1}{2}{\theta^\prime}^2+\frac{1}{2}Q\sin^2\theta+\left(1-Q\right)\cos\theta \label{reduced1}
\end{equation}

\noindent and the equation of motion is

\begin{align}
	\frac{d}{d\tau}\frac{\partial \mathbf{L_1}}{\partial\theta^\prime}=\theta^{\prime\prime}=\frac{\partial \mathbf{L_1}}{\partial\theta}&=Q\sin\theta\cos\theta-(1-Q)\sin\theta \\
	& \nonumber \\
	\therefore\theta^{\prime\prime}&=-\left[1-Q-Q\cos\theta\right]\sin\theta \label{rpp.eom}
\end{align}

\subsection{Equilibria and their Stability}

The equilibria of the system are given by

\begin{align}
	1-Q-Q\cos\theta&=0 \text{ ,}\quad\textbf{OR}\;\; \theta=0,\pm\pi \nonumber \\
	\cos\theta&=\frac{1-Q}{Q} \nonumber \\
	\theta=\pm\sec^{-1}&\left(\frac{Q}{1-Q}\right)  \text{ for } Q\geq\frac{1}{2}
\end{align}

\begin{figure}[!ht]
	\includegraphics[scale=0.7]{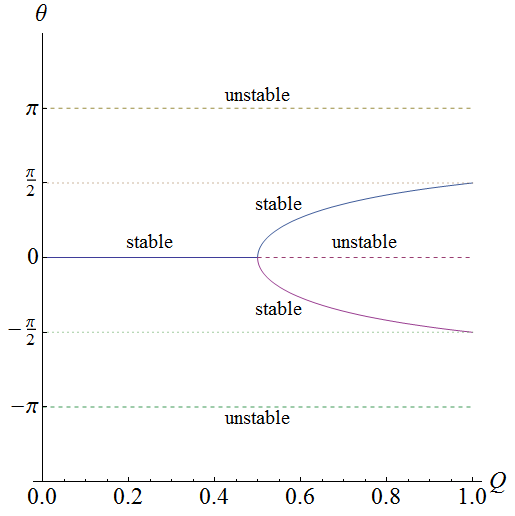}\label{rpp.biff}
	\caption{The equilibria of a rotating physical pendulum, and the bifurcation at parameter value $Q=1/2$.}
\end{figure}

Fig. \ref{rpp.biff} shows the stability of the equilibria, as is easily confirmed by the potential function $V(\theta)$, and it's second derivative

\begin{align}
	V(\theta)=&-\frac{Q}{2}\sin^2\theta-(1-Q)\cos\theta \nonumber \\
	\nicefrac{d^2V}{d\theta^2}=&-2Q\cos^2\theta+(1-Q)\cos\theta+Q\,.
\end{align}

The linearization of (\ref{rpp.eom}) about any equilibrium $\theta=\theta^*$ is

\begin{equation}
	(\theta-\theta^*)^{\prime\prime}=-\left(\left.\frac{d^2V}{d\theta^2}\right\rvert_{\theta=\theta^*}\right)(\theta-\theta^*)
\end{equation}

\noindent so that near equilibrium, solutions are an oscillatory normal mode $\theta-\theta^*\sim e^{\pm i\Omega\tau}$ if the second derivative of $V$ is positive, and $\Omega=\sqrt{\nicefrac{d^2V}{d\theta^2}}$.  Alternately, if the second derivative of $V$ is negative, then the equilibrium is a saddle and has two exponential modes; one growing, one decaying:  $\theta-\theta^*\sim e^{\pm\,\alpha\tau}$, where $\alpha=\sqrt{-\nicefrac{d^2V}{d\theta^2}}$.  In general we refer to $\sqrt{\abs{\nicefrac{d^2V}{d\theta^2}}}$ as the normal mode rate (NMR).

To better understand the physics of these rates, consider again the definition (\ref{def.omeg}) of $\omega$, and let $\omega_g^2=\nicefrac{K}{(M\ell^2+I_N)}$, $\omega_r^2=\nicefrac{(M\ell^2+I_\perp-I_P)\omega_a^2}{(M\ell^2+I_N)}$ define gravitational and rotational characteristic frequencies, respectively, such that $\omega^2=\omega_r^2+\omega_g^2$. We then have that

\begin{align}
	Q\omega^2=\omega_r^2\;\;, \quad (1-Q)\omega^2=\omega_g^2
\end{align}

\noindent and $Q=\omega_r^2/(\omega_r^2+\omega_g^2)$.  The bifurcation $Q=1/2$ and conditions such as $Q>1/2$ correspond to $\omega_r=\omega_g$ and $\omega_r>\omega_g$, respectively.  We use these relations to evaluate the final column of Table~\ref{nmf}, where we also use the alternate parameter $q=Q/(1-Q)=\nicefrac{\omega_r^2}{\omega_g^2}=\bar{E}/K$.

\begin{table}[!ht]
	\centering
	\caption{Classification and Normal mode rates of the Equilibria of a Rotating Physical Pendulum }
	\begin{ruledtabular}
		\begin{tabular}{l  p{5.5cm} p{5.8cm}}
			Equilibrium &  dimensionless NMR $\sqrt{\abs{\nicefrac{d^2V}{d\theta^2}}}$ & physical $\Omega\omega$, $\alpha\omega$ \\
			\hline
			$\theta=0$                & $\Omega=\sqrt{1-2Q}$ for $0\leq Q<1/2$, $\alpha=\sqrt{2Q-1}$ for $1/2<Q\leq1$ & $\Omega\omega=\sqrt{\omega_g^2-\omega_r^2}\text{ for }\omega_r<\omega_g,\quad\quad\;$ $\alpha\omega=\sqrt{\omega_r^2-\omega_g^2}\text{ for }\omega_r>\omega_g$  \\
			\hline
			$\theta=\pm\pi$           & $\alpha=1$ & $\alpha\omega=\sqrt{\omega_r^2+\omega_g^2}$ \\
			\hline
			$\theta=\pm\sec^{-1}(q)$ & $\Omega=\sqrt{2-Q^{-1}}$ for $1/2<Q\leq1$ & $\begin{array}{l}
				\Omega\omega=\omega_r\sqrt{1-\omega_g^4/\omega_r^4}\text{ for }\omega_r>\omega_g\\
				\hphantom{\Omega\omega}=\omega_r\sqrt{1-q^{-2}}\text{ for }q>1
			\end{array}$ \\
		\end{tabular}
	\end{ruledtabular}                        
	\label{nmf.tab}
\end{table}

\section{Derivation of the RDP Lagrangian}\label{Ap.RDP.deriv}

\begin{figure}
	\includegraphics[scale=.25]{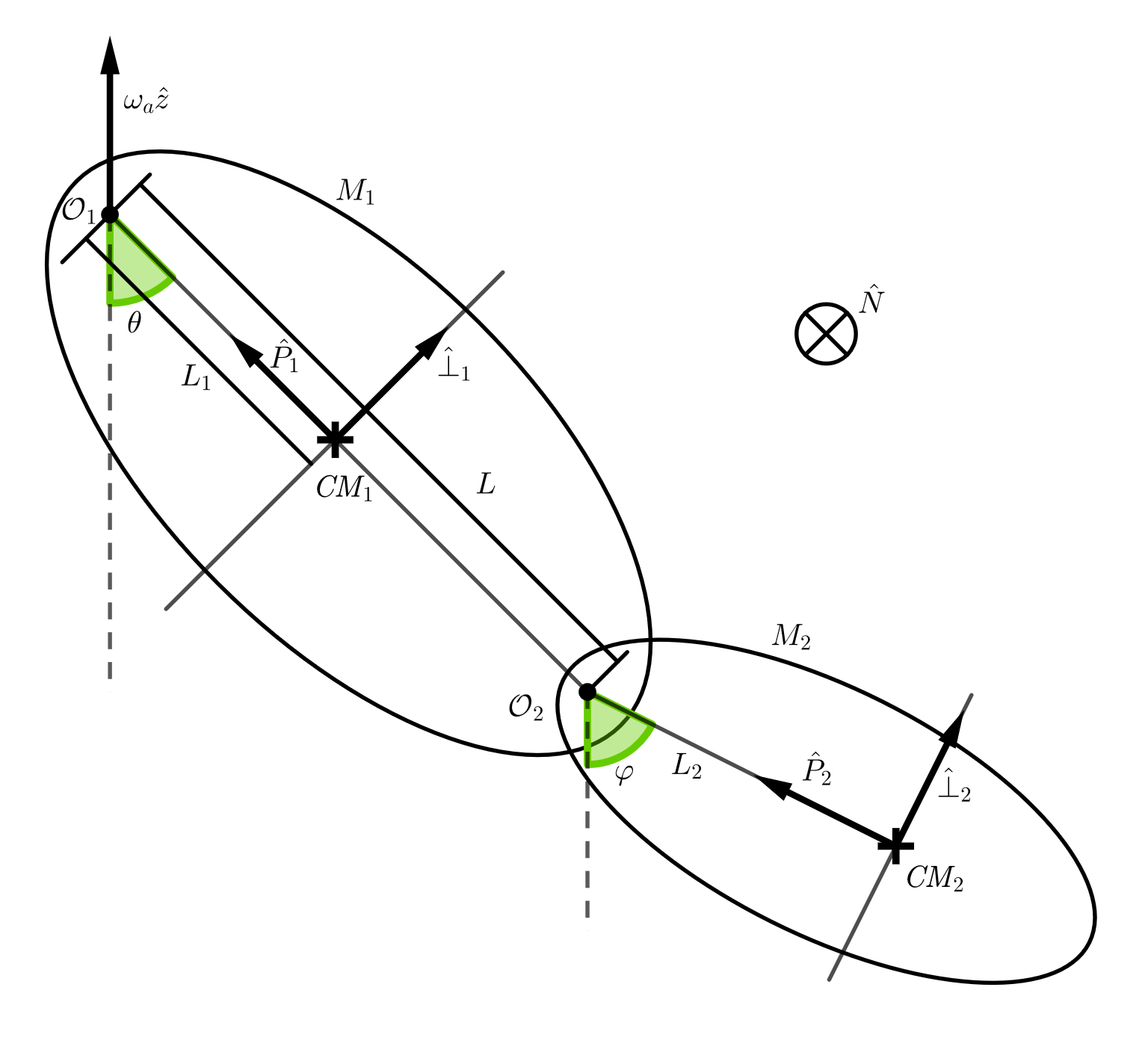}
	\caption{Schematic showing the construction and configuration of a rotating double pendulum from two physical pendula.}
	\label{RDPconfig}
\end{figure}

In this Appendix we present the details of the derivation of the RDP Lagrangian.  We begin with the rotating physical pendulum as described in Appendix~\ref{rpp}, and attach a second rigid body, as shown in Fig.~\ref{RDPconfig}.  The upper pendulum has physical parameters (that is, with units) as follows: $M_1,L_1,I^{(1)}_{P,\perp,N}$.  Pivot the second pendulum at a point $\mathcal{O}_2$ on the $P_1$ axis, so that the pivots $\mathcal{O}_1, \mathcal{O}_2$ and the center of mass $\mathit{CM}_{\!1}$ of the first pendulum are collinear.  The pivots are a distance $L$ apart.  The second pendulum is similarly aligned via its principal axes (its pendulum axis $P_2$ is an angle $\varphi$ away from vertical), and has the following parameters: mass $M_2$, $L_2$ is the distance from $\mathcal{O}_2$ to the center of mass $\mathit{CM}_{\!2}$ of the second pendulum, and principal moments $I^{(2)}_{P,\perp,N}$. The Lagrangian is

\begin{equation}
	\mathcal{L}_2=\mathcal{L}_1+\frac{1}{2}\vec{\omega}_2\cdot[I^{(2)}]\vec{\omega}_2+\frac{1}{2}M_2v_{\mathit{CM}_{\!2}}^2+M_2\,g\left(L\cos\theta+L_2\cos\varphi\right)\,.
\end{equation}

\noindent The angular velocity $\vec{\omega}_2$ is

\begin{equation*}
	\vec{\omega}_2=\omega_a\hat{z}-\dot{\varphi}\hat{N}=\left(\omega_a\cos\varphi,\omega_a\sin\varphi,-\dot{\varphi}\right)\,,
\end{equation*}

\noindent from which the rotational kinetic energy is

\begin{equation}
	\frac{1}{2}\vec{\omega}_2\cdot[I^{(2)}]\vec{\omega}_2=\frac{1}{2}I^{(2)}_P\omega_a^2\left(1-\sin^2\varphi\right)+\frac{1}{2}I^{(2)}_\perp\omega_a^2\sin^2\varphi+\frac{1}{2}I^{(2)}_N\dot{\varphi}^2\,.
\end{equation}

The velocity of $\mathit{CM}_{\!2}$ has a component in the vertical plane, and an $\hat{N}$-component due to rotation $\omega_a$.  The component in the plane is the velocity of $\mathcal{O}_2$ plus the velocity of $\mathit{CM}_{\!2}$ relative to $\mathcal{O}_2$.  The velocity due to rotation about the vertical involves the horizontal distance from $\mathcal{O}_1$ to $\mathit{CM}_{\!2}$. This gives

\begin{equation}
	\vec{v}_{\mathit{CM}_{\!2}}=\underbrace{L\dot{\theta}\hat{\perp}_1}_{\vec{v}_1}+\underbrace{L_2\dot{\varphi}\hat{\perp}_2}_{\vec{v}_2}+\underbrace{(L\sin\theta+L_2\sin\varphi)\omega_a\hat{N}}_{\vec{v}_3}\,.
\end{equation}

\noindent Observe that $\hat{\perp}_1\cdot\hat{\perp}_2=\cos\left(\varphi-\theta\right)$ and $\hat{\perp}_i\cdot\hat{N}=0$, so

\begin{align}
	v_{\mathit{CM}_{\!2}}^2=&v_1^2+v_2^2+v_3^2+2\vec{v}_1\cdot\vec{v}_2 \nonumber \\
	=&L^2\dot{\theta}^2+L_2^2\dot{\varphi}^2+\left(L\sin\theta+L_2\sin\varphi\right)^2\omega_a^2 \nonumber\\
	&+2L_2L\cos\!\left(\varphi-\theta\right)\dot{\theta}\,\dot{\varphi}\,.
\end{align}

\noindent As with the rotating physical pendulum in Appendix~\ref{rpp}, we non-dimensionalize time $t=\omega\tau$, so that the kinetic terms change by $\dot{\theta}\rightarrow\omega\theta^\prime$ and $\dot{\varphi}\rightarrow\omega\varphi^\prime$, $\omega$ to be determined.  The Lagrangian of the rotating double pendulum is

\begin{align}
	\mathcal{L}_2=&\frac{1}{2}A\,{\theta^\prime}^2+B\cos\!\left(\varphi-\theta\right)\theta^\prime\,\varphi^\prime+\frac{1}{2}C\,{\varphi^\prime}^2 \nonumber\\
	&+\frac{1}{2}\bar{A}\sin^2\theta+\bar{B}\sin\theta\sin\varphi+\frac{1}{2}\bar{C}\sin^2\varphi \label{rdp.lag.ap}\\
	&+K_1\cos\theta+K_2\cos\varphi \nonumber
\end{align}

\noindent where the coefficients are

\begin{align}
	A&=\left(M_1 L_1^2+M_2L^2+I^{(1)}_N\right)\omega^2 \nonumber\\
	B&=M_2 L_2 L \omega^2 \nonumber\\
	C&=\left(M_2 L_2^2+I^{(2)}_N\right)\omega^2 \nonumber\\
	\bar{A}&=\left(M_1L_1^2+M_2L^2+I^{(1)}_\perp-I^{(1)}_P\right)\omega_a^2\quad  \label{coeffs.ap}\\
	\bar{B}&=M_2 L_2 L \omega_a^2 \nonumber\\
	\bar{C}&=\left(M_2 L_2^2+I^{(2)}_\perp-I^{(2)}_P\right)\omega_a^2 \nonumber\\
	K_1&=\left(M_1L_1+M_2L\right)g \nonumber\\
	K_2&=M_2L_2g\,. \nonumber
\end{align}

\noindent The effective potential energy is

\begin{align}
	U_\text{eff}=-&\frac{1}{2}\left(\bar{A}\sin^2\theta+2\bar{B}\sin\theta\sin\varphi+\bar{C}\sin^2\varphi\right) \nonumber\\
	&-K_1\cos\theta-K_2\cos\varphi\,.
\end{align}

The coefficients are constrained once a characteristic frequency $\omega$ is chosen.  Following the method in Appendix~\ref{rpp}, we proceed by establishing a connection between the Lagrangians $\mathcal{L}_2$ and $\mathcal{L}_1$, namely:  $\mathcal{L}_2$ has the same form as $\mathcal{L}_1$ if the pendulum angles $\theta,\varphi$ are forced to coincide.

\begin{align}
	\mathcal{L}_1=&\left.\mathcal{L}_2\right\vert_{\varphi\rightarrow\theta} \nonumber\\
	=&\frac{1}{2}\left(A+2B+C\right){\theta^\prime}^2+\frac{Q}{2}\left(\bar{A}+2\bar{B}+\bar{C}\right)\sin^2\theta+(K_1+K_2)\cos\theta
\end{align}

Thus we identify $E=A+2B+C$, $\bar{E}=\bar{A}+2\bar{B}+\bar{C}$ and $K=K_1+K_2$, which gives

\begin{align}
	\omega^2=\frac{\bar{E}+K}{M_1 L_1^2+M_2(L+L_2)^2+I^{(1)}_N+I^{(2)}_N}
\end{align}

\noindent and the coefficients of the RDP Lagrangian (\ref{rdp.lag}) are constrained by

\begin{equation}
	A+2B+C=\bar{A}+2\bar{B}+\bar{C}+K_1+K_2\,.
\end{equation}

\section{Full Outputs of Polynomial Systems}\label{outputs}

\subsection{The General Case}\label{outputs1}
The basic set of equations is the system for equilibrium; other systems are built from this by adding further equations.  The polynomial system for equilibrium of the RDP consists of

\begin{align}
	\frac{\partial V}{\partial\theta}=0&\,\,\implies\,\,\text{\lstinline|- 2*s1 - 2*chi*s1 + c1*d*qq*s*s1 + 2*c1*qq*s2 - c1*d*qq*s2|} \\
	\frac{\partial V}{\partial\varphi}=0&\,\,\implies\,\,\text{\lstinline|- 2*s2 + 2*chi*s2 + 2*c2*qq*s1 - c2*d*qq*s1 + 2*c2*d*qq*s2|}\\
	&\qquad\qquad\text{\lstinline|- c2*d*qq*s*s2|} \nonumber \\
	\intertext{and we also include}
	&\text{\lstinline|c1^2 + s1^2 - 1|} \\
	&\text{\lstinline|c2^2 + s2^2 - 1|}
\end{align}

\noindent where \texttt{chi} stands for $\chi$, \texttt{d} for $1+\delta$, \texttt{s} for $1+\sigma$, and \texttt{qq} for $q$.  The bifurcation system is the above equations together with the Hessian determinant of the potential

\begin{flalign}
	&\quad\det\left(H[V](\theta,\varphi)\right)=0\quad\implies &&
\end{flalign}

\vspace{-12pt}

\begin{lstlisting}
	4*c1*c2 - 4*c1*c2*chi^2 - 4*c1*c2^2*d*qq - 4*c1*c2^2*chi*d*qq - 4*c1^2*c2^2*qq^2 + 4*c1^2*c2^2*d*qq^2 - c1^2*c2^2*d^2*qq^2 - 2*c1^2*c2*d*qq*s + 2*c1*c2^2*d*qq*s + 2*c1^2*c2*chi*d*qq*s + 2*c1*c2^2*chi*d*qq*s + 2*c1^2*c2^2*d^2*qq^2*s - c1^2*c2^2*d^2*qq^2*s^2 + 2*c2*d*qq*s*s1^2 - 2*c2*chi*d*qq*s*s1^2 - 2*c2^2*d^2*qq^2*s*s1^2 + c2^2*d^2*qq^2*s^2*s1^2 + 4*c1*qq*s1*s2 + 4*c2*qq*s1*s2 + 4*c1*chi*qq*s1*s2 - 4*c2*chi*qq*s1*s2 - 2*c1*d*qq*s1*s2 - 2*c2*d*qq*s1*s2 - 2*c1*chi*d*qq*s1*s2 + 2*c2*chi*d*qq*s1*s2 - 4*c2^2*d*qq^2*s1*s2 + 2*c2^2*d^2*qq^2*s1*s2 - 2*c1^2*d*qq^2*s*s1*s2 + 2*c2^2*d*qq^2*s*s1*s2 + c1^2*d^2*qq^2*s*s1*s2 - c2^2*d^2*qq^2*s*s1*s2 + 2*d*qq^2*s*s1^3*s2 - d^2*qq^2*s*s1^3*s2 + 4*c1*d*qq*s2^2 + 4*c1*chi*d*qq*s2^2 - 2*c1*d*qq*s*s2^2 - 2*c1*chi*d*qq*s*s2^2 - 2*c1^2*d^2*qq^2*s*s2^2 + c1^2*d^2*qq^2*s^2*s2^2 + 4*qq^2*s1^2*s2^2 - 4*d*qq^2*s1^2*s2^2 + d^2*qq^2*s1^2*s2^2 + 2*d^2*qq^2*s*s1^2*s2^2 - d^2*qq^2*s^2*s1^2*s2^2 + 4*d*qq^2*s1*s2^3 - 2*d^2*qq^2*s1*s2^3 - 2*d*qq^2*s*s1*s2^3 + d^2*qq^2*s*s1*s2^3
\end{lstlisting}

\subsection{PMMR Equations}\label{output2}

Elimination $\chi$ in the above system by (\ref{chiPMMR}) results in the following system

\begin{align}
	&\text{\lstinline|- 4*d*s*s1 + 2*c1*d*qq*s*s1 - c1*d^2*qq*s*s1 + c1*d^2*qq*s^2*s1 + 4*c1*qq*s2|} \nonumber \\
	&\text{\lstinline|- 4*c1*d*qq*s2 + c1*d^2*qq*s2 + 2*c1*d*qq*s*s2 - c1*d^2*qq*s*s2,|} \\
	\nonumber \\
	&\text{\lstinline|4*c2*qq*s1 - 4*c2*d*qq*s1 + c2*d^2*qq*s1 + 2*c2*d*qq*s*s1 - c2*d^2*qq*s*s1|}\nonumber \\
	&\text{\lstinline|- 8*s2 + 4*d*s2 + 4*c2*d*qq*s2 - 2*c2*d^2*qq*s2 - 2*c2*d*qq*s*s2 |} \\
	&\text{\lstinline|+ 3*c2*d^2*qq*s*s2 - c2*d^2*qq*s^2*s2,|}\nonumber\\
	\nonumber \\
	&\text{\lstinline|c1^2 + s1^2 - 1,|} \\
	&\text{\lstinline|c2^2 + s2^2 - 1,|}
\end{align}


\begin{figure}[!h]
	\includegraphics[scale=0.85]{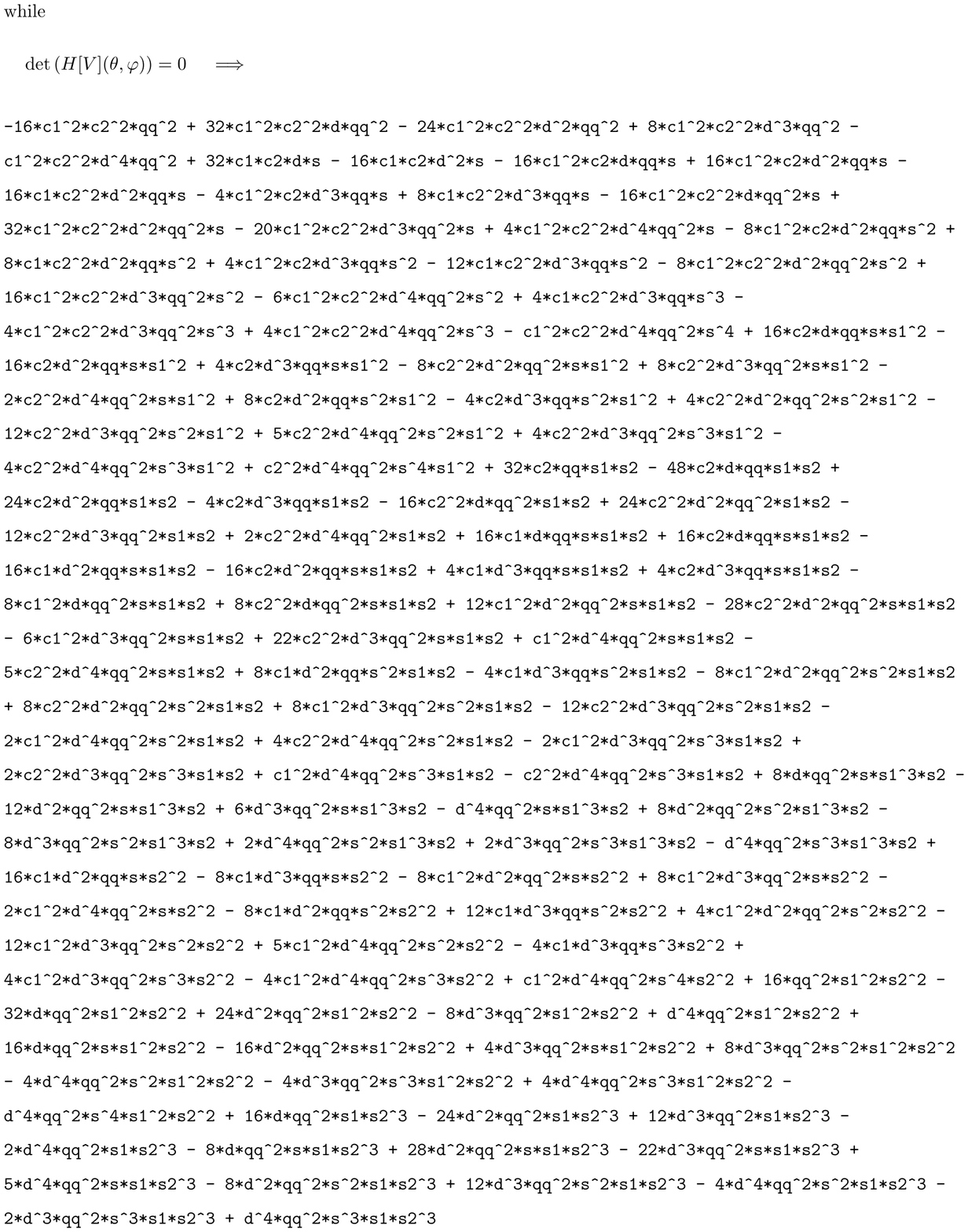}
\end{figure}

\end{document}